\documentclass[12pt,leqno]{amsart}
\newtheorem{theorem}{Theorem}[section]
\newtheorem{lemma}[theorem]{Lemma}
\newtheorem{corollary}[theorem]{Corollary}

\begin{document}

\title{Commutants of Bol Loops of Odd Order}
\author{Michael~K.~Kinyon}
\address{Department of Mathematical Sciences \\
Indiana University South Bend \\
South Bend, IN 46634 USA}
\email{mkinyon@iusb.edu}
\urladdr{http://www.iusb.edu/\symbol{126}mkinyon}
\author{J.~D.~Phillips}
\address{Department of Mathematics \& Computer Science \\
Wabash College \\
Crawfordsville, IN 47933 USA}
\email{phillipj@wabash.edu}
\urladdr{http://www.wabash.edu/depart/math/faculty.htm{\#}Phillips}

\begin{abstract}
In this note we show that the commutant of a Bol loop
of odd order is a subloop.
\end{abstract}

\maketitle
\thispagestyle{empty}

\section{Introduction}
\label{sec:intro}

A \textit{loop} $(L,\cdot)$ is a set $L$ with a binary operation
$\cdot : L\times L\to L$ such that (i) for given $a,b\in L$, the equations
$a\cdot x = b$ and $y\cdot a = b$ have unique solutions $x,y\in L$,
and (ii) there exists a neutral element $1\in L$ satisfying $1\cdot x
= x\cdot 1 = x$ for all $x\in L$. Basic references for loop theory are
\cite{bruck}, \cite{pflug}. We will use the usual juxtaposition
conventions to avoid excessive parenthesization, e.g., $ab\cdot c
= (a\cdot b)\cdot c$. The \textit{commutant} of a loop $L$
is the set
\[
C(L) = \{ a\in L : ax = xa\ \ \ \forall x\in L\} .
\]
The \textit{center} of $L$ is the set of all $a\in C(L)$ such that
$a\cdot xy = ax\cdot y$, $x\cdot ay = xa\cdot y$, and $xy\cdot a = x\cdot ya$
for all $x,y\in L$. The center is a normal subloop. For some varieties of loops,
such as groups, the commutant and center coincide. For other varieties, the
commutant is larger than the center, but is still ``well-behaved" in the sense
that it is a normal subloop. However, the commutant of an arbitrary loop need
not be a subloop at all, and even when it is, it need not be normal.

The commutant is also known in the literature by other names such as
``Moufang center", ``commutative center", or ``centrum". Since this object
is not, in general, central in the sense of universal algebra, we prefer a term that
does not suggest otherwise. Thus we have borrowed the term ``commutant",
which is used for a similar concept in other fields.

A loop is said to be a (left) \textit{Bol loop} if it satisfies the identity
$x(y\cdot xz) = (x\cdot yx)z$ for all $x,y,z$. A right Bol loop is similarly
defined, and a loop which is both a left and right Bol loop is a \textit{Moufang loop}.
(This is one of many equivalent definitions; see \cite{bruck}, \cite{pflug}, \cite{robinson}.)
In this paper, all Bol loops will be left Bol loops. The commutant of a Moufang loop is
a subloop, but it is an open problem to characterize precisely those Moufang
loops for which the commutant is normal \cite{doro}. On the other hand, the commutant
of an arbitrary Bol loop need not be a subloop; see the web page \cite{moorhouse}. (Note
that \cite{moorhouse} lists right Bol loops.) Our main result is the following.

\begin{theorem}
\label{thm:main}
If every element of the commutant of a Bol loop has finite odd order 
then the commutant is a subloop.
\end{theorem}

Glauberman observed that a finite Bol loop has odd order
if and only if every element has odd order \cite{glauberman}. Thus we have the following.

\begin{corollary}
\label{coro:main}
The commutant of a finite Bol loop of odd order is a subloop.
\end{corollary}

Our investigations were aided by the automated reasoning tool OTTER,
developed by McCune \cite{otter}. We thank Ken Kunen for introducing
us to OTTER.

\section{Proofs}
\label{sec:proofs}
Let $L$ be a Bol loop. We recall that $L$ has both the \textit{left inverse property}
(LIP) and the \textit{left alternative property} (LAP) \cite{pflug} \cite{robinson}.
These are given, respectively, by the equations $x^{-1}\cdot xy = x\cdot x^{-1}y = y$
and $x\cdot xy = x^2 y$, where $x^{-1}x = xx^{-1} = 1$. We shall make use of
these without comment in what follows.

\begin{lemma}
\label{lem:intersect}
Let $L$ be a Bol loop. If $a,b$ are in $C(L)$,  then so are $a^2$, $b^{-1}$ and $a^2 b$.
\end{lemma}

\begin{proof}
Fix $x\in L$. We have $a (a^2 x) = (a^2 x) a = (a\cdot ax) a = (a\cdot xa) a = a\cdot x a^2$.
Cancelling, we have $a^2 x = x a^2$, that is, $a^2 \in C(L)$. Next, we have
$b^{-1}x = b^{-1}(x\cdot b^{-1}b) = (b^{-1}\cdot xb^{-1})b = b(b^{-1}\cdot xb^{-1}) = xb^{-1}$,
so that $b^{-1}\in C(L)$. Finally,
$
a^2 b\cdot x = (a\cdot ba)x = a(b\cdot ax) = a(b\cdot xa) = a(b\cdot x(b\cdot b^{-1}a))
= a((b\cdot xb)\cdot ab^{-1}) =  (a\cdot (b\cdot xb)a)b^{-1} = (a^2 \cdot (b\cdot xb)) b^{-1}
= ((b\cdot xb)\cdot a^2) b^{-1} = (b(x\cdot ba^2))b^{-1} = b^{-1}\cdot b(x\cdot ba^2) =
x\cdot a^2 b .
$
Therefore $a^2 b \in C(L)$.
\end{proof}

For an element $a$ of a loop $L$, let $\langle a\rangle$ denote the subloop
generated by $a$. Bol loops are \textit{power-associative} \cite{pflug} \cite{robinson},
that is, if $x^0 := 1$, $x^{n+1} := xx^n$, and $x^{-n-1} := x^{-1} x^{-n}$ for
$n \geq 0$, then $x^m x^n = x^{m+n}$ for all integers $m, n$.

\begin{lemma}
\label{lem:powers}
Let $L$ be a Bol loop. For each $a\in C(L)$, $\langle a\rangle \subseteq C(L)$.
\end{lemma}

\begin{proof}
Since $L$ is power-associative, this follows from Lemma \ref{lem:intersect}
and induction.
\end{proof}

We now can prove Theorem \ref{thm:main}.
Given $a, b\in C(L)$, let $c\in L$ be the unique element such that
$c^2 = a$. Since $a$ has odd order, Lemma \ref{lem:powers} implies that
$c\in C(L)$. By Lemma \ref{lem:intersect}, $ab = c^2 b\in C(L)$.
Since $C(L)$ is closed under products and inverses, we may
apply the \textit{left inverse property} $a^{-1}\cdot ab = b$,
which holds in Bol loops  \cite{pflug} \cite{robinson}, and the
commutativity of $C(L)$ to conclude that $C(L)$ is a subloop.

We conclude with three remarks. First, the previously cited result
of Glauberman actually applies to any subset of a Bol loop containing
the neutral element and closed under taking inverses and under the
operation $(x,y)\mapsto x\cdot yx$ \cite{fkp}. So in fact, we have
the following extension of Corollary \ref{coro:main}: \textit{If the
commutant of a Bol loop has finite odd order, then the commutant
is a subloop.}

Second, an element $a$ of an arbitrary loop $L$ is called a \textit{Bol element}
if $a (x\cdot ay) = (a \cdot xa)y$ for all $x,y\in L$. The set $B(L)$ of all Bol elements
of $L$ need not be a subloop. The proofs herein carry over nearly verbatim to obtain
the following extension of Theorem \ref{thm:main}: \textit{If every element of}
$B(L)\cap C(L)$ \textit{has finite odd order, then} $B(L)\cap C(L)$ \textit{is a subloop}.

Finally, we have been unable to find an example of a (necessarily infinite) uniquely $2$-divisible
Bol loop (i.e., a loop in which the mapping $x\mapsto x^2$ is a bijection) with a commutant
element whose unique square root is not a commutant element.  However, we conjecture
that such loops do indeed exist.

\end{document}